
\input amssym
\input epsf
\magnification\magstep1
\font\deffont=cmssi10

\def\de#1{{\deffont #1}}
\def\complexes{{\Bbb C}}
\def\reals{{\Bbb R}}
\def\integers{{\Bbb Z}}
\def\naturals{{\Bbb N}}

\def\bn{{\bf n}}
\def\bX{X}
\def\bY{{Y}}
\def\curl{\hbox{\rm curl}}
\def\area{\hbox{\rm area}}
\def\divv{\hbox{\rm div}}
\def\grad{\nabla}

\def\raw{\rightarrow}

\def\dalpha{{\dot \alpha}}
\def\section#1#2{\medskip\noindent{\bf \S #1 #2.}}
\def\frac#1#2{{#1 \over #2}}
\def\thm#1{\medskip\noindent{\bf #1}}
\def\thmp#1#2{\medskip\noindent{\bf #1} (#2)}

\def\laplacian{\Delta}
\def\br{{\bf r}}
\def\balpha{(\alpha_i)}
\def\phit{{\phi_t}}
\def\bibitem#1{\item{[#1]}}
\def\QED{  \rlap{$\sqcup$}$\sqcap$ \medskip}

\centerline{\bf Dynamics of two-dimensional time-periodic Euler fluid flows}
\medskip
\centerline{Philip Boyland}
\centerline{Department of Mathematics}
\centerline{University of Florida}
\centerline{Gainesville, FL 32605-8105}
\centerline{boyland@math.ufl.edu}
\bigskip
{\narrower

\noindent{\bf Abstract:}  This paper investigates
the dynamics of time-periodic Euler flows in multi-connected, planar 
fluid regions which are ``stirred'' by the moving boundaries. The classical 
Helmholtz theorem on the transport of vorticity implies that if
the initial vorticity of such a flow is generic among real-valued functions 
in the $C^k$-topology ($k \geq 2$) or is $C^\omega$ and nonconstant,
then the flow has zero topological entropy.  On the other hand, it
is shown that for constant initial vorticity there are stirring protocols
which always yield time-periodic Euler flows with positive entropy.
These protocols are those that generate flow maps in pseudoAnosov
isotopy classes. These classes are a basic ingredient of the 
Thurston-Nielsen theory and a further application of that theory 
shows that pseudoAnosov stirring protocols with generic initial vorticity 
always yield solutions to Euler's equations for which
the sup norm of the gradient of the vorticity grows exponentially in
time.  In particular, such Euler flows are never time-periodic.

}

\bigskip
\section{1}{Introduction}
The motion of fluids provides intuition and
basic terminology for Dynamical Systems theory, and 
the ideas and results of low dimensional dynamics continue to have many
important applications in Fluid Mechanics. 
In this paper we explore a particular case of a general question 
connecting Fluid Mechanics and Dynamical Systems:
\thm{Question 1:}  {\it
Are fluid flows typical dynamical systems, or do they have
special, distinguishing dynamical features?
}
\medskip
This question is obviously
vague and can be made precise in a number of ways depending
on the dimension, fluid model and time dependence of the
velocity fields. In this paper we study the dynamics of
the Poincar\'e maps of time-periodic Euler flows
(Question 2 in \S 3).

There is a great deal known Question 1
for steady (i.e. time-independent)  Euler flows 
in two and three dimensions. Much of this work is summarized
in the book of Arnol'd and Khesin ([AK]). In the interval
since its publication significant progress
has been made by Etnyre and Ghrist on steady 3D Euler flows (see [GK]
for a summary).  In particular, they have shown that 
every nonsingular $C^\omega$-Euler flow on the 
three sphere has a periodic orbit that is unknotted ([EG]).  This is a feature
of Euler flows which is not shared by every volume preserving
flow and so provides a distinguishing dynamical feature for
the class of 3D steady Euler flows. 

It is a common heuristic that
the dynamics of surface diffeomorphisms share many of the
features of 3D flows, and so it is natural from the
Dynamical System's standpoint to study the time $T$-maps of
$T$-periodic 2D Euler flows. The lower dimension results in
topological simplification, but the time periodicity
adds a new feature. In addition, we shall primarily be concerned
with the situation in which the region occupied by the fluid changes 
during its evolution as the fluid is ``stirred'' by various 
protocols.  Thus the time $T$-flow maps can be in different isotopy classes, 
and the insights provided by the Thurston-Nielsen
theory of these classes will be central to our analysis. 

This paper views Fluid Mechanics from the Dynamical System's
prospective. From this admittedly skewed point of view,
Fluid Mechanics is the study of families of diffeomorphisms (the flow maps) 
which have the special property  that their
vector fields satisfy equations such as Navier-Stokes,
Euler or Stokes. These equations are,
of course, derived from physical principles, and the resulting systems
are studied because they have been found to agree very well with the behavior
of real fluids. Question 1 asks 
whether these physically based conditions restrict the possible
dynamics a fluid can manifest as compared to a general
dynamical system.

For example, in many physical situations the amount of expansion
or compression of a fluid during its evolution is very
small and can be neglected. Thus 
fluid maps are often assumed to be area or volume preserving. 
As is well known, the dynamics of
area-preserving systems differ from general dynamical
systems in many respects, for example, they have no attractors.
Thus in formulating the corresponding  precise version of Question 1 
in two dimensions 
one should ask whether the dynamics of fluids are 
distinguished from those of a general area-preserving system. 

The version of Question 1 studied here assumes that 
the fluid motion satisfies the Euler equations. These equations are
based on the assumption of a fluid with no viscosity and
thus there is no energy dissipation.  While we restrict attention here 
to $C^k$-classical solutions, lesser regularity and
various kinds of weak solutions have been an object
of much study and are of central importance.
Classical results about Euler fluid motions and their
consequences for the dynamics are given in \S 3 and \S 4.
Theorem 2 states a fundamental result of 
Helmholtz and Kelvin from the mid $19^{th}$ century. It  
says that Euler fluid motions are characterized by their 
preservation of vorticity (or curl)  and circulation integrals. These
distinguishing features are metric dependent, but nonetheless provide a
key to essential topological properties of Euler fluid motions.
Proposition 4 shows that as a consequence of the 
preservation of vorticity, 
a time-periodic Euler fluid motion whose initial
vorticity is a typical $C^r$-function or else is $C^\omega$ and
nonconstant must have zero topological entropy (the $C^\omega$ case
is equivalent to Remark 4 in [BS]). 

Amongst the Euler flows with non-generic vorticity those with
zero vorticity are of particular interest.
These systems are much studied as a consequence of their 
mathematical tractability and because it is often argued that 
a fluid at rest has zero vorticity and so if the vorticity is
conserved, it will be zero for all time.  Using
standard potential theory, Proposition 3 shows that systems
with constant vorticity and periodic stirring protocols always give
rise to periodic Euler fluid motions. Since conditions that ensure
periodic solutions to the Euler equations are very rare, this
result is useful in guaranteeing the existence of at least one interesting
class of time $T$-Euler flow maps. 

Using these results in the real 
analytic case yields a dichotomy that is 
somewhat similar to that of Arnol'd for steady 3D Euler flows (see II \S 1 
in [AK]).  For 2D time-periodic Euler flows, if the vorticity is nonconstant,
then the dynamics are ``integrable'' and the entropy is zero. One
may have chaotic dynamics in the constant vorticity case, but  only
if the fluid is stirred

The next section, \S 5, contains results with the opposite conclusion
from that of Proposition 3.
As a consequence of the Thurston-Nielsen theory of surface automorphisms,
Theorem 7 shows that  for generic initial vorticity there are large classes of
periodic stirring protocols which never give rise to periodic
Euler fluid motions. Further, 
these stirring protocols result in exponential
growth of the sup norm of the gradient of the vorticity.
As a consequence of the preservation of vorticity 
one expects this general type of behavior in 
any chaotic Euler flow with generic smooth vorticity.
The attractive feature of Theorem 7 is that the topology of the
stirrer motion allows one to get concrete results on the exponential
growth.

It is anticipated that most of this paper's readers will be 
familiar with Dynamical Systems and not Fluid Mechanics,
so we have made some attempt to include basic Fluid Mechanics 
from a Dynamical Systems point of view.  A second reason
for including statements and proofs of classic results
is that the case under consideration here, multi-connected
fluid regions with moving boundaries, is not usually discussed in
the standard texts.

There are many first class books on mathematical Fluid Mechanics.
The books [AK], [Ba], [C], [G], [MB], [MP], and [Se] 
provide a sample with a variety of emphases.
For a survey of Fluid Mechanics and mathematical structures
with a point of view similar to this paper, see [Bd2].
The motion of planar fluids stirred by topologically complex
protocols has been studied in [BAS1], [FCB], and [V].

\section{2}{Fluid motions and Hamiltonian systems}
In this paper a \de{fluid motion}, $(M_t, \phit)$, consists of 
a smooth family of smooth planar fluid regions, $M_t$,
with the Euclidean metric and a smooth one-parameter family of
diffeomorphisms, $\phi_t :M_0 \rightarrow M_t$, with $\phi_0 = id$. 
The diffeomorphisms should be thought
of as describing the evolution of fluid particles: the particle at the point
$z \in M_0$ at time $0$ is at the point $\phi_t(z) \in M_t$
at time $t$. 

In our fluid region families the outer boundary of $M_t$ will always be 
a fixed, smooth simple closed curve $C_0$.
In addition, $m$  disks, each of radius $\epsilon$, 
(the \de{stirrers}) are excluded
from the fluid region, with the case $m=0$ being allowed. 
The evolution of the fluid region is determined 
by the rigid motion of the excluded disks.
The stirrers always remain
circles of radius $\epsilon$ and their centers move along 
paths given by $\alpha_i: [0, \infty)  \raw D$, for $i = 1, \; \dots\; ,
m$.  It is assumed that no two disks collide and
no disk collides with the outer boundary. 
Because of the slip boundary conditions for Euler fluid motions, 
rotations of the boundary do not affect the fluid motion and so
will be ignored.  The collection of paths $ \balpha =
(\alpha_1,\; \dots\;, \alpha_m)$ is called a \de{stirring protocol}.
Since the outer boundary remains fixed, specifying
a stirring protocol specifies  the fluid regions $M_t$, and vice versa.
The inner boundary circles at time 
$t$ are denoted $C_{it}$ and have velocity $\dalpha_i$. 
To simplify notation 
we let $\alpha_0$ be a constant function, so the
velocity of $C_{0}$ is identically zero: $\dalpha_0 \equiv 0$.

The \de{velocity} or \de{vector field} $\bX$ of the fluid motion 
at the point $\phi_t(z)$ at time $t$ is defined by
$$ \bX(\phi_t(z), t) := \frac{\partial \phi_t}{\partial t}(z).\eqno{(1)}$$
Note that in contrast to what is usual in Dynamical Systems,  
the vector field $X$ may, and usually does, depend on time. Thus 
$\phi_t$ is a flow in the usual Dynamical System
sense if and only $\bX$ is time-independent.  In that case the 
fluid motion is called  \de{steady}. A family of smooth vector field
vector fields $X(\cdot, t)$ on $M_t$ can always be integrated to give a family 
of diffeomorphisms $\phit$. The fact that the fluid at
time $t$ is contained solely in $M_t$, i.e. 
the image of $\phit$ is $M_t$, is equivalent to the boundary
conditions 
$$\bX \cdot \bn_i =  \dalpha_i \cdot \bn_i, \eqno{(2)}$$ 
where $\bn_i$ is the unit normal to the boundary circle $C_{it}$.
Without comment in the sequel we will go back and forth between a family of
velocity fields satisfying these conditions and the corresponding
diffeomorphisms and will use both $(M_t, \phit)$ and $(M_t, X)$ to denote
the fluid motion.

At this point the notion of a fluid motion is very
general and it is not required to satisfy any particular equation.
The first restriction on a fluid motion comes from 
assuming that the fluid neither expands nor contracts 
during its evolution.  A fluid motion is called \de{incompressible} 
if each diffeomorphism $\phi_t$ in the fluid motion is 
area-preserving, or equivalently,  if $\divv(\bX) \equiv 0$. 

In two dimensions area preservation is closely related 
to being a Hamiltonian system.
Recall that the family $\phi_t$ is called \de{Hamiltonian}
if there exist a family of real-valued functions $H_t$, so that the
the velocity field $X$ satisfies $\bX(\cdot, t) = J\;\grad H_t $ 
for each $t$, where
$J = \pmatrix{0 & 1 \cr -1 & 0 \cr}$. In this case each diffeomorphism
$\phit$ is also called \de{Hamiltonian}.
In many cases, a system being Hamiltonian forces additional
dynamical properties beyond those of area preservation.
(There is a vast literature arising from the Arnol'd conjecture
and subsequent developments. See, for example, [MS] and [FH].)
However, the distinction between area preservation and 
Hamiltonian on surfaces is associated with the flux across generators of
homology that are not associated with boundary curves. Thus the
distinction disappears in genus zero. This is the content of 
first part of the next lemma. We include it here as it introduces a
standard construction of the stream function which we shall
need later, and because the motion of the boundary curves introduce
an element not contained in the usual results.
The second part of the lemma concerns a related question: Is every
area preserving diffeomorphism $f : M_0 \raw M_0$ (which may permute the
inner boundaries) Hamiltonian?

\thm{Lemma 1:} {\it

(a) The fluid motion $(M_t, \phit)$ 
 is Hamiltonian if and only if it is incompressible.

(b) If $f:M_0 \raw M_0$ is an
area-preserving diffeomorphism, then there is a Hamiltonian 
fluid motion  $(M_t, \phit)$  with $f = \phi_1$.
}

\thm{Proof:} For part (a), first note that
Hamiltonian systems are always area preserving.
For the converse, 
at each fixed time, on the boundary component 
$C_{it}$ the vector $\dalpha_i$ is a constant and so by direct calculation,  
$$
0 = \oint_{C_{it}} \dalpha_i \cdot d\bn_i. \eqno{(3)}
$$
Thus if for each $t$ we let $\bY = -J\, \bX$,  then 
$\curl(\bY) =\divv(\bX) =  0$ and using (2) and (3),
$0 = \oint_{C_{it}} \bY \cdot d\br $.
So for each $t$, $\bY$ is a curl-free field
which has zero circulation around each boundary curve. 
Thus we may fix $z_0$ and unambiguously define
$$\Psi(z, t) = \int_{z_0}^z \bY \cdot d\br \eqno{(4)}$$ 
which will satisfy $\bY = \grad \Psi$.
Therefore,  $\bX = J\, \grad \Psi$, which is to say that
$\Psi$ is the Hamiltonian that generates $\phi_t$.  

For part (b), first note that since $M_0$ is a disk with holes
removed it is standard  that we may find a $1$-periodic
protocol $\balpha$ with corresponding regions $M_t$ and a 
fluid motion ${\hat\psi_t}:M_0 \raw M_1$ with 
${\hat \psi_1}$ isotopic to $f$ on $M_0 = M_1$ ([Bi]). 
In addition, since for each $t$ the region $M_t$ has the same area as $M_0$, 
using a theorem of Moser ([M]) we may find a diffeomorphism $g_t:M_t \raw
M_t$, isotopic to the identity, with 
$g_{t*} ({\hat\psi_{t*}}\, \nu_0) = \nu_t$, where 
$\nu_t$ is the Euclidean area form on $M_t$. Note that Moser's
theorem is stated for closed manifolds, but the proof also works for
manifolds with boundary and yields a family 
$g_t$ which is smooth in $t$.  If we let $\psi_t = g_t \circ {\hat \psi_t}$,
then $(M_t, \psi_t)$ is an area-preserving fluid motion.

By construction, $\psi_1^{-1} \circ f$ is area preserving
and isotopic to the identity
on $M_0$. By a similar argument using Moser's theorem (cf. remark 1.4.C in
[P]) we may find a family of area preserving diffeomorphisms $h_t:M_0 \raw
M_0$ with $h_0 = id$ and $h_1 = \psi_1^{-1} \circ f$. The incompressible
and thus Hamiltonian fluid motion required for (b) is 
$\phit = \psi_t \circ h_t$.  \QED

In Fluid Mechanics the Hamiltonian $\Psi$ is called the 
\de{stream function}.  The role of Lemma 1 in this paper 
is to show that in formalizing 
Question 1 for our fluid regions it suffices to consider
area-preserving systems; we do not need to say that 
we are looking for Euler models of a general time-periodic Hamiltonian 
system.

\section{3}{Euler fluid motions and diffeomorphisms}
In its simplest interpretation, the Euler equation is Newton's second
law, $F = ma$, applied to each fluid particle with a force
resulting from the gradient of the pressure. As is common,
we assume a uniform fluid density of one.
What we call an Euler fluid motion here
should properly be called an incompressible, constant density
Euler flow. We omit the adjectives as understood and use the terminology
``fluid motion'' to avoid confusion with the Dynamical Systems notion
of flow.  Euler flow is also sometimes called perfect or ideal fluid flow.
From this point onward all fluid motions are required to have 
velocity fields which are $C^k$ for $k \geq 3$.

\thm{The Euler equations for the fluid motion:} 
{\it An incompressible fluid motion $(M_t, \phit)$ 
is called Euler 
if there is a smooth family of smooth
functions $p_t: M_t \rightarrow \reals$ (the pressure)
so that $p_t$ and $\phi_t$ satisfy
$$
\frac{\partial^2 \phit}{\partial t^2}(z)  = -\grad p_t(\phi_t(z)).
$$
}

It is more conventional to express the Euler equations in term
of the velocity field.  The acceleration term becomes 
the \de{material derivative} of the vector field defined as  
$$
\frac{D\bX}{Dt} := \frac{\partial^2 \phit}{\partial t^2} 
= \frac{\partial (\bX\circ\phit) }{\partial t} 
= \frac{\partial \bX}{\partial t} + \nabla_{\bX} \bX,
$$
where $\nabla_\bX \bX$ is  the directional or covariant  
derivative of $\bX$ in the direction of $\bX$.

\thm{The Euler equations for velocity fields:} {\it 
A fluid motion $(M_t, X)$ is a 
solution to the incompressible Euler equation
if there is a smooth
family of smooth functions $p_t: M_t \rightarrow \reals$ (the pressure)
so that $p_t$ and $\bX$  satisfy
$$\eqalign{   
\frac{D\bX}{Dt} &=  - \grad p_t \hbox{} \cr 
\divv(\bX) &= 0 \cr
\bX \cdot \bn_i & =  \dalpha_i \cdot \bn_i\ \hbox{\rm on the boundary.}\cr}
$$
}

It is easy to check that a 
fluid motion is incompressible Euler if and only if it
velocity field is a solution to the incompressible Euler equations.
Since the fluid regions, $M_t$, are determined by the stirring
protocol, $\balpha$, we shall also say that the Euler fluid motion
$(M_t, X)$ is an Euler solution which is \de{compatible}
with the given stirring protocol.

The basic questions of existence,  uniqueness, and regularity
for solutions
to the Euler equations has been much studied (see [C], [MB], [MP], 
or [Kh] for expositions of various results). For the case
of interest here, multiply connected planar regions with moving
boundary,  existence and uniqueness of a global weak and classical 
solutions were obtained by Kozono [Ko] and He and Hsiao [HH].

The class of dynamical systems to be studied here 
are the time $T$-diffeomorphisms of
time-periodic Euler fluid motions. A stirring protocol
is said to be \de{$T$-periodic} if after time $T$, the set of
stirrers return to their initial position,
and then the stirring process repeats. 
Note that since the stirrers are indistinguishable, we only 
require the set of stirrers to come back to itself, not each
individual stirrer. 
In terms of the paths $\alpha_i$, this definition  says that 
$\alpha_i(t + T) = \alpha_{\sigma(i)}(t)$ for some
permutation $\sigma$ of $ \{ 1, 2, \; \dots\;, m\} $. 

\thm{Definition of Euler Diffeomorphisms:} {\it 
A diffeomorphism $f : M_0 \rightarrow M_0$ is called
an \de{Euler diffeomorphism}, if $f = \phi_T$ for  
an incompressible Euler fluid motion $(M_t, \phi_t)$ 
whose velocity field $X(z, t)$ is
$T$-periodic and compatible with a $T$-periodic stirring protocol.
}
\medskip

Note that the definition requires not just 
a periodic protocol, but the compatible Euler 
velocity field must be periodic as well. 
The usual definition of a Hamiltonian diffeomorphism requires
that $f = \phi_1$ for a Hamiltonian family $\phit$. This in,
in fact, equivalent to requiring $f = \phi_1$ for a family
$\phit$ that is generated by a family of Hamiltonians $H_t$
that are $1$-periodic, $H_{t+1} = H_t$, and so the
corresponding vector fields $X(z, t)$ are also $1$-periodic
(see, for example, page 37 in [P]).
On the other hand, as a consequence of Theorem 7 there are 
diffeomorphisms, $f$, that are the time $1$-map of an
Euler fluid motion, but that Euler fluid motion cannot have
a periodic velocity field. The general question of which 
area preserving diffeomorphisms are the time $1$-maps of perhaps 
aperiodic Euler fluid motions is very difficult, but
see Brenier [Br] and Shnirelman [Sh] for  reports on significant progress.
 The requirement
that an Euler diffeomorphism be generated by a periodic
velocity field is included in our definition to insure that the
iterates of $f$ describe the dynamics of the fluid motion.

We can now give a precise statement of the version of Question 1
of interest here.

\thm{Question 2:} {\it 
Given an area preserving diffeomorphism $g:M_0 \raw M_0$, is there
an Euler diffeomorphism $f$ that is topologically conjugate 
to $g$?}
\medskip

As is usual in Dynamical Systems, we only require
$f$ and $g$ to have the same dynamics up to topological change
of coordinates. This is essential because being an Euler fluid motion
is not even preserved under smooth changes of coordinates because
the acceleration, or equivalently, the covariant derivative, depends
on the metric.

The next theorem collects classical results 
that are basic to understanding the dynamics of Euler fluid motions. 
The \de{vorticity} of a velocity field is it
curl and is denoted $\omega_t(z) := \curl(\bX(z, t))$. In
two dimensions the vorticity  is a real-valued function and 
$\omega_t = -\laplacian \Psi$, where 
$\Psi$ is the stream function of $\bX(z, t)$.
Recall that the push forward of a scalar field ($0$-form) $s$ under
a diffeomorphism $f$ is $f_* s  = s \circ f^{-1}$

\thmp{Theorem 2:}{Helmholtz-Kelvin} {\it
An incompressible fluid motion $(M_t, \phi_t)$ 
with velocity field $\bX$ and vorticity $\omega_t$ is Euler
if and only if its vorticity is passively transported, 
$$
\phi_{t\ast}\,\omega_0 = \omega_t, \eqno{(5)}
$$
and circulations around all smooth
simple closed curves $C$ are preserved under the flow, 
$$
\frac{d}{dt}\oint_{\phi_t(C)} \bX \cdot d\br = 0.\eqno{(6)}
$$
}
\medskip

\thm{Proof.} 
First observe that standard two-dimensional vector identities yield
$$
\frac{D\bX}{Dt} = 
\frac{\partial \bX}{\partial t} + \grad( \frac{\;\|\bX \|^2}{2}) -
\omega J \bX 
$$
Since $\curl (J \bX) = \divv( \bX)  = 0$, taking the curl yields, 
$$
\curl\left( \frac{D\bX}{Dt}\right)  = 
\frac{\partial \omega}{\partial t} +
\grad \omega \cdot \bX = 
\frac{\partial \omega(\phi_t(z), t)}{\partial t}. \eqno{(7)}
$$
Now if $\bX$ is Euler,  then $\curl (\frac{D\bX}{Dt})  = 
- \curl( \grad p)  = 0$, and 
so by (7), $\omega$ is constant on orbits, which is equivalent to (5).

The transport theorem for simple closed curves $C$ is
$$
\frac{\partial}{\partial t} \oint_{\phi_t(C)} \bX \cdot d\br =
\oint_{\phi_t(C)} \frac{D\bX}{Dt} \cdot d\br. \eqno{(8)}$$
So if  $\bX$ is Euler,  
$$
\frac{\partial}{\partial t} \oint_{\phi_t(C)} \bX \cdot d\br 
= -\oint_{\phi_t(C)} \grad p \cdot d\br = 0,
$$
and circulations are preserved.

For the converse, first note that if vorticity is transported, 
then (7) yields
$\curl (\frac{D\bX}{Dt})  = 0.$
In addition, if circulation integral
are conserved then (8) yields
$$0 = \oint_{\phi_t(C)} \frac{D\bX}{Dt} \cdot d\br.$$
Thus $\frac{D\bX}{Dt}$ is a curl-free field whose circulation around
each boundary component is zero. Thus as in (4), for each $t$ 
there exists a smooth function $-p_t$ with 
$\frac{D\bX}{Dt} = -\grad p_t$, and further, since all  the
data is varying smoothly with $t$, so does $p_t$.
\QED

By Green's theorem, if $C$ bounds a disk in $M_0$, then 
(5) implies (6),  but for multi-connected regions (6) is stronger.
In particular, the requirement that circulation integrals be preserved
is needed for multi-connected regions. This condition is
often not included in standard texts in the vorticity form of the Euler 
equations because the regions under consideration
are usually simply connected or else solutions with discontinuous
pressure are allowed.

\section{4}{Constant and generic initial vorticity}
In this section we take the first steps on
Question 2 with a pair of results that give dynamical information
about Euler fluid motions based on the nature of
the initial vorticity distributions. Since vorticity distribution
coupled with the circulations determine the velocity field
and the vorticity is passively transported by the Helmholtz-Kelvin
Theorem, it is not surprising
that the initial vorticity distribution very much influence
the dynamics. 

In view of the definition of Euler diffeomorphism
we restrict attention to periodic stirring protocols.
As noted above, in general,  a periodic stirring protocol does not
suffice to insure that a compatible Euler solution is periodic.
However, as a consequence of elementary potential theory 
it does suffice if the initial vorticity is constant. This is
content of the first part of Proposition 3.
The second part of the proposition concerns the
situation with stationary boundaries and is a special case of
a more general well-known classical result (see, for example, 
Proposition 8.2.2 in [AMR], or page 33 in [MP]).

\thm{Proposition 3:} 
{\it 
Assume that $\balpha$ is a $T-$periodic stirring protocol with
corresponding fluid regions $M_t$. Given 
a real number $\Omega$ and a vector 
$(\Gamma_0, \Gamma_1,  \; \dots\;, \Gamma_m)
\in \reals^{m+1}$ 
with $\sum \Gamma_i = \Omega\;  \area(M_0) $,  there exists 
a unique $T$-periodic incompressible Euler fluid motion $(M_t, X)$ with
$\omega_0 \equiv \Omega$ and $\oint_{C_{i0}} \bX \cdot d\br = \Gamma_i$,
for $i = 0, \dots, m$.
In particular, if the inner boundaries are stationary,
then the solution $\bX$ is steady (time-independent) 
and so $h_{top}(\phi_t) = 0$ for all $t$.
}

\thm{Proof:}  Fix $t$.
By standard arguments (eg. a simple modification of the proof
of Theorem 2.2 in [MP]), given the planar 
region with smooth boundaries $M_t$, the
circulations $\Gamma_i$ with $\sum \Gamma_i  = \Omega\; \area(M_t) $, 
and the vectors $\dalpha_i$, 
there is a unique $C^\omega$-function $\Psi_t$ with
$\laplacian \Psi_t = -\Omega$, 
$(J\,\grad \Psi) \cdot \bn_i = \dalpha_i \cdot \bn_i$ on each boundary
component $C_{it}$, and
$\oint_{C_{it}} (J\;\grad \Psi) \cdot d\br = \Gamma_i$ for
all $i$. Note that the condition on the sum of 
the circulations is necessary by Green's  theorem.

Now define $\bX = J\; \grad \Psi_t$, and let $\phi_t$ be it's
fluid motion. By construction, $\phi_t$ transports vorticity and
preserves circulation integrals and so by Theorem 2
is an Euler fluid motion. Since the $M_t$ are a
$T$-periodic family, $\bX$ is $T$ periodic. Uniqueness follow
from uniqueness of each $\Psi_t$ and Theorem 2.

If the boundaries of the fluid region
do not move,  then $M_t = M_0$ for all $t$,
and so $\Psi_t = \Psi_0$ for all $t$, and so $\bX$ is time-independent.
Thus $\phi_t$ is a two-dimensional flow
in usual Dynamical Systems sense, and so has zero topological entropy ([Y]).
\QED

The next result concerns the dynamics  when the initial
vorticity is typical amongst all
functions in the $C^k$ category, $2 \leq k \leq \infty$. 
The basic idea, as was observed by Brown and Samelson [BS], is that as a 
consequence of the Kelvin-Helmholtz Theorem, an
Euler diffeomorphism preserves the level sets of its initial vorticity
distribution. Since smooth, one-dimensional manifolds cannot
support chaotic dynamics, if the level sets are smooth manifolds, the
diffeomorphism has zero topological entropy. 
The precise version of this idea 
we give makes use of a frequently invoked theorem of Katok.

\thm{Proposition 4:} {\it If $f$ is an
Euler diffeomorphism whose velocity field $\bX$ has initial vorticity 
$\omega_0$ which has finitely many critical points, 
then $h_{top}(f) = 0$. 
}

\thm{Proof.}
By Theorem 1, and the fact that $\bX$ is $T$-periodic,
$\phi_{T*}\, \omega_0 = \omega_T = \omega_0$. 
Now assume to the contrary that
$h_{top}(f) > 0$.  As a consequence of a theorem of Katok [K2], 
since $f$ is smooth, it must have a transverse homoclinic point to a 
hyperbolic periodic point $p$. Now since $\omega_0 = \omega_0\circ f$, and 
$\omega_0$ is continuous, it must be constant on the closure of
the union of  the stable and unstable manifolds of the orbit of $p$. 
And since these intersect transversally at the homoclinic orbit and
$\omega_0$ is smooth, every point on the homoclinic orbit must
be a critical point for $\omega_0$, and so there are infinitely
many of them, contrary to assumption.
\QED

The class of functions allowed for $\omega_0$ in the proposition contains
the Morse functions, and so contains an open dense set
in  $C^k(M_0, \reals)$, $2 \leq k \leq \infty$. Thus it represents the typical
case in the $C^k$-topology.  It important to note what the proposition
does not say. It does not show that amongst the initial vorticities
that give rise to periodic Euler solutions, the typical case
has zero topological entropy. Such a stronger result would require
understanding conditions on the initial vorticity that insure or
even allow periodicity of the Euler velocity field (cf. Question 3 in
\S 6).  Also note that any condition on a smooth
invariant function which guarantees zero entropy can be
substituted for the assumption of finitely many critical points.
Finding such conditions is an interesting general Dynamical
Systems question.

A nonconstant real analytic function has finitely
many critical points in any compact set and so we have a 
corollary that via Katok's theorem is equivalent to Remark
4 in [BS].

\thmp{Corollary 5:}{Brown and Samelson}  {\it If $f$ is an
Euler diffeomorphism whose velocity field is $C^\omega$ and
has nonconstant initial vorticity, then $h_{top}(f) = 0$.}

\section{5}{PseudoAnosov stirring protocols}
In this section we use the fact that the topological character of
a $T$-periodic stirring protocol determines the isotopy class
of the time $T$-map of any compatible fluid motion, and this information
can be used with the Thurston-Nielsen theory to obtain classes of
periodic stirring protocols for which any compatible Euler solution
is never periodic.

We begin by quickly reviewing standard material 
about  braids, isotopy classes and Thurston-Nielsen theory. 
For more information on braids and isotopy classes see [Bi] or [BL].
For general material on the Thurston Nielsen theory see [T], [FLP], and
[CB].  For a survey of dynamical applications of the TN theory see [Bd1], 
and for fluid mechanical applications see [BAS1] and [BAS2].

For definiteness we fix the geometry of the fluid
regions.  In this section $M_0$ is always the unit disk 
with three small round disks of
radius $\epsilon$ centered at $p_1 = -1/2$, $p_2 = 0$, and $p_3 = 1/2$ 
removed:
$$
M_0 := \{ z \in \complexes : 
|z| \leq 1, \ \hbox{\rm and}\ |z - p_i| \geq \epsilon \ \hbox{\rm for} \
i = 1, 2, 3 \}.
$$
We shall also assume that any stirring protocol is $1$-periodic. Such a 
protocol can be associated with braid on three strands. 
 This standard construction proceeds by considering the three-dimensional
traces of the paths defined by $(\alpha_i(t), t)$ for $0 \leq t \leq 1$ and
$i = 1, 2, 3$.
This defines a physical braid on three strands. By examining the crossing
of the strands projected onto a plane, we obtain a braid word in
the braid group on three strings, $B_3$. 
Note that changing the projection plane changes the braid word
by a conjugacy, and so the construction actually only defines a 
conjugacy class in $B_3$, which is known as the \de{braid type} of
the protocol. The distinction between the braid and its conjugacy
class is usually not of consequence here, so we will just refer to
the ``braid'' of the protocol.

The braid in turn determines the isotopy class in the mapping
class group of $M_0$. The Thurston-Nielsen (TN) theory contain
a classification of surface isotopy classes into three
types: finite order, pseudoAnosov (pA), and reducible. A
stirring protocol is said to have one of these types if
its corresponding isotopy class does. Since
we restrict attention here to $M_0$, we can use a simplified
version of the theory which depends on the fact that $M_0$ 
is a two fold cover of the torus with $4$ disks removed ([Bi], [BW], 
[K1], see \S 1.8 in [Bd1] for an exposition).

The TN-type of an isotopy class represented by an element
of $B_3$ can be computed directly using the homomorphism 
$\chi: B_3 \raw SL(2, \integers)$ with
$$
\chi(\sigma_1) = \pmatrix{1 & 1 \cr 0 & 1} \ \ \ 
\chi(\sigma_2) = \pmatrix{1 & 0 \cr -1 & 1} 
$$
Note that this is just the Burau representation with 
the substitution $t = -1$, and so can be interpreted as
the action on homology in the two-fold cover, the torus.   
The result we need is that the braid word $\beta \in B_3$ represent
a pA class if and only if the largest eigenvalue $\lambda'$ of the
matrix $\chi(\beta)$ is real with magnitude larger than one, 
or equivalently, if $|\hbox{\rm trace}(\chi(\beta))|
> 2$. The number $\lambda = |\lambda'|$ is called
the \de{expansion constant} of the pA class.

PseudoAnosov isotopy classes have a number of
special properties.  Of central importance here is the
fact that a pA isotopy class 
always induces exponential word growth on $\pi_1$. More
precisely, let $F_3 \cong \pi_1(M_0)$ be the free group on
three generators and $\rho:F^3 \raw F^3$ be induced by a
pA isotopy class, i.e.  $\rho = f_*:\pi_1(M_0) 
\raw \pi_1(M_0)$ for any (and thus every) homeomorphism $f$ in the class.
If for $w \in F_3$, $\ell(w)$ denotes the
number of letters in the reduced word for $w$,
then for every nontrivial $w \in F_3$, 
$$\lim_{n\raw\infty} \ell(\rho^n(w))^{1\over n} = \lambda, \eqno{(9)}
$$
where to be explicit, $\rho^n(w)$ denotes repeated composition
of $\rho$ applied to $w$, not multiplication in the group.

One consequence of this induced growth on $\pi_1$ is that any homeomorphism
in the class has topological entropy greater than
$\log(\lambda)$. A second consequence concerns essential
curves and arcs. An \de{essential simple closed curve}
is one that is neither contractible nor boundary parallel. 
An arc with its endpoints on the boundary is \de{essential} if
cannot be contracted to a point by a homotopy that keeps its endpoints
on the boundary. As a consequence of (9), under a homeomorphism in
a pA class, no iterate of an essential arc or curve is homotopic 
to itself. Here as in the rest of the paper
all homotopies are required to keep the boundary fixed set-wise.

A third consequence of (9) is the exponential growth
of the length of essential curves and arcs under iteration. 
If $\gamma$ is a smooth essential arc or closed curve and $f_1, f_2,
\;\dots \;$ is a sequence of 
diffeomorphisms all in the same pA isotopy class with expansion 
constant $\lambda$, then there is a positive
constant $k$ so that 
$$\ell((f_1\circ f_2\;\dots\circ f_n)(\gamma)) \geq k \lambda^n \eqno{(10)}$$
for all $n > 0$, where $\ell(\cdot)$ now denotes the Euclidean length. 
The equivalence of (9) and (10) follows 
from the fact that the word length of an element
in $\pi_1$ gives a uniform bound on  the displacement caused
by its corresponding deck transformation in the universal
cover (see [Mi1] and [FLP], expos\'e 10, \S II). This result is
usually stated for the repeated application of a single 
diffeomorphism, but since the composition above induces the same
action on $\pi_1$ as a repeated composition, the 
proof is the same.  For what follows it
is important to note that the constant $k$ depends only
on the homotopy class of the essential arc or closed curve.

In order to get the estimate in Theorem 7 we require (generic) hypothesis
on the initial vorticity. 
Recall that $\mu: M_0 \raw \reals$ is a \de{Morse function} if
its Hessian is nonsingular at  all critical points. 
Morse functions are open and dense in $C^k(M_0, \reals)$, $2\leq k \leq\infty$.
(see [Mi2], [H], or [Ma] for more details on Morse Theory). 
A {\it connected component} of a level set $\mu^{-1}(r)$ is
called a \de{critical set} if it contains a critical point and
a \de{regular} set if it does not.
Regular sets are either simple closed curves or arcs with both endpoint
on the boundary of $M_0$.

Informally, an \de{essential regular strip} is a collection
of parallel regular arcs. More precisely, an essential regular strip
is defined to be a nontrivial interval
of $\mu$ values $I$  and a collection of 
essential regular arcs $K_r$ for each $r\in I$ 
so that each $K_r$ is a connected component
of the level set $\mu^{-1}(r)$ and $S := \cup  K_r$
is compact and connected. An \de{essential regular annulus} 
is defined similarly using regular closed curves. 
If $I = [c, d]$, the boundary of an essential regular annulus consists
of the two regular closed curves $K_{c}$ and $K_{d}$. 
The boundary of an essential regular strip is the union
of four arcs: 
the two segments in the boundary of $M_0$ that consists of
the endpoints of the regular arcs in the strip  and  
the two regular arcs $K_{c}$ and $K_{d}$.

Note that if $K$ is an essential regular arc (resp. closed curve), 
then it always has
a neighborhood (one-sided, if $K$ is tangent to the boundary)
that is an essential regular strip (annulus).
All regular sets in an essential regular strip or annulus are 
in the same homotopy class.  It is worth noting that 
not every Morse function on $M_0$ has essential regular sets, see
Figure 1a.  

\midinsert
\def\epsfsize#1#2{.8\hsize}
\centerline{\epsfbox{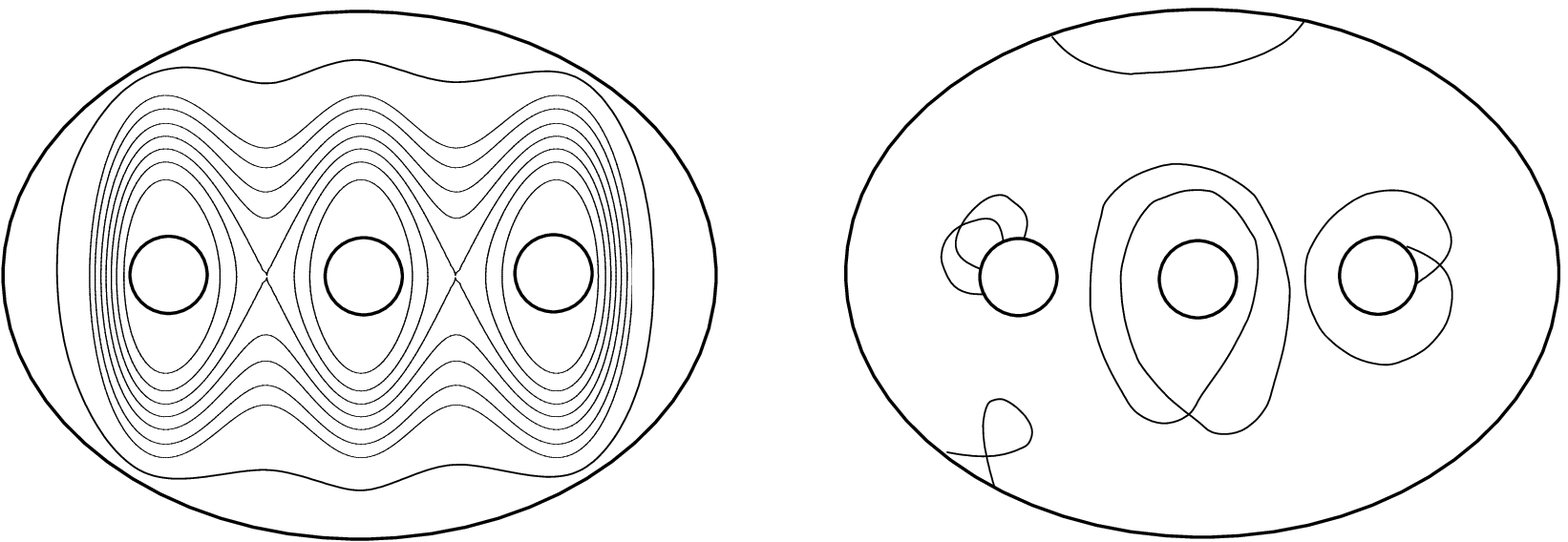}}
\medskip
{\leftskip=50pt\rightskip=50pt\parindent=0pt

{\bf Figure 1a:}  Sample level sets of a Morse function with no 
essential regular set. \hfill\break
{\bf Figure 1b:}  Components of level sets from  the
proof of Lemma 6

}

\endinsert

\thm{Lemma 6:} {\it  Let $M_0$ be the disk with three holes as defined above.
For all $2 \leq k \leq \infty$ there is a dense, open set 
$G\subset C^k ( M_0, \reals)$ such that each $\mu\in G$ has an 
essential regular set.
}
 
\thm{Proof:}
Let the set $G$ consists of Morse functions $\mu \in C^k ( M_0, \reals) $
for which 

{\narrower

\item{(1)} Each critical value comes from exactly one critical point, i.e. 
if $c_1$ and $c_2$ are distinct critical points, then $\mu(c_1) \not = 
\mu(c_2)$. 

\item{(2)} No critical point is on the boundary and no critical set
is tangent to the boundary. 

\item{(3)} Restricted to each boundary circle $\mu$ is a Morse function
with one-dimensional domain.

}

Standard arguments show that $G$ is open dense
in the Morse functions, and thus is open dense in 
$ C^k ( M, \reals).$
By condition (1) each critical set of a $\mu\in G$ 
contains exactly one critical point and thus there are three
possibilities for the topology of a critical set:
an `$x$' with the endpoints of all
its arms on the boundary, an `$\alpha$' with the endpoint
of its two arms on the boundary, and an `$\infty$'. 
By condition (2)  any critical set that
intersects the boundary must be an '$\alpha$' or '$x$' with 
its endpoints on the boundary, and 
by condition (3) a given level set can only
intersect the boundary in a finite set of points.

Now assume contrary to the conclusion of the lemma 
that $\mu\in G$ has no essential regular set. 
We claim that this implies that for each boundary circle, $C_i$, there
is a regular closed curve, $B_i$,  that is homotopically parallel to $C_i$.
The proof proceeds by drawing a number of conclusions about 
the possible level sets of $\mu$ and their configurations all
based on the same general argument:
if the assertion is not true, then near to the level set under 
consideration there
would be another component of a level set that is an essential
arc or closed curve, contrary to the assumption on $\mu$. For example, 
the endpoints of the arms of a critical set of
type '$\alpha$' must be on the same
boundary component, for if they were not, we could move a little off
the essential path in the $\alpha$ from endpoint to endpoint and
find an essential arc for $\mu$. Other assertions that follow from
the same argument are:
A critical set of type '$\alpha$'  must either be contractible into a boundary
circle or else contain its legs inside its closed loop and
that closed loop is homotopically parallel to the boundary
circle (see Figure 1b).  A critical set of type 
'$x$' must have all its endpoints on the same boundary circle
and must be contractible into that circle (see Figure 1b).
A regular arc must have both of its endpoints on the same boundary
curve, be contractible into that boundary circle,
and not be tangent to any other boundary circle. 
A regular closed curve 
can be tangent to one and only one boundary circle, but then it 
must be contractible into that boundary circle 
or else be homotopically parallel to it.

Now for $i = 0,\; \dots\; 3$, let 
$$F_i = \{ K : K \ \hbox{\rm is a component of a level set and}\ 
K\cap C_i \not = \emptyset \}.$$ 
By construction, each $F_i$ is connected and using the 
observations of the previous paragraph,  for $i \not = j$,
$F_i \cap F_j  = \emptyset$. In addition, 
$U = (\cup F_i)^c$ is open because a component of
a level set not intersecting a boundary curve implies that
all level sets near it also do not intersect the boundary curve.
Thus each $F_i$ is a compact, connected set that is disjoint from
the other $F_j$. This implies that the topological frontier of each
$F_i$  consists of a finite number of components of
level sets which are disjoint from the other $F_i$.
Looking at the list of possible elements of an $F_i$ we
see that this frontier can only be a closed curve tangent and
parallel to $C_i$  or else a critical set of '$\alpha$'-type
with its arms on $C_i$ and its loop parallel to $C_i$ as on the right 
of Figure 1b. In either case, just outside the frontier there is a regular
closed curve, $B_i$, that is parallel to $C_i$, proving the claim above.

If we let $M'$ be the multi-connected region whose boundary is  the union 
of the $B_i$, then $M'$ is
homeomorphic to $M_0$ and has regular closed curves  for
its boundary. As a consequence,  all the components of level sets of $\mu$ in
$M'$ are either regular closed curves or else  critical sets of
type '$\infty$'.  For $i = 0,\; \dots\; 3$, let 
$E_i$ be the closure of the set of regular closed curves $K$ which
are homotopically parallel to the boundary circle $C_i$ and focus
on a fixed $i$. 
The frontier of $F_i$ gives one component of the frontier of $E_i$. 
The rest of the frontier of $E_i$ must consist of components of critical 
sets of type '$\infty$' and one of the loops of the $\infty$ must
be parallel to $C_i$.  If the other loop of the $\infty$ 
bounds a disk or is parallel to $C_i$, then nearby the $\infty$ there is
a regular closed curve, parallel to $C_i$, contradicting the fact that
the $\infty$ was on the boundary of $E_i$. On the other hand,
if the other loop of the $\infty$ encloses a boundary circle other than 
$C_i$, then nearby there
would be an essential regular closed curve, contrary to the assumption
on $\mu$. Thus in either case we have a contradiction to the
assumption that $\mu$ has no regular sets.
\QED

The next result gives conditions under which a periodic
stirring protocol gives rise to Euler fluid motions which
are never periodic. On one hand, by Proposition 4, if the
initial vorticity has finitely many critical points then
any Euler diffeomorphism has zero topological entropy. On
the other hand, any diffeomorphism in a pA isotopy class must 
have positive entropy.  Thus any Euler fluid motion whose vorticity
has finitely many critical points and is compatible
with a pA stirring protocol cannot have a periodic 
velocity field. With the additional generic condition
on the initial vorticity given in Lemma 6, one gets more data about this
non-periodicity in the form of an estimate of the
growth rate of the gradient of the vorticity.

\thmp{Theorem 7:}{PseudoAnosov protocols} {\it
Assume that $\balpha$ is a stirring protocol of pseudoAnosov type
with expansion constant $\lambda$, and $\bX$ is a compatible Euler
solution with initial vorticity $\omega_0$. If $\omega_0$ has finitely
many critical points, then $\bX$ is not periodic. If in
addition, $\omega_0\in G$, the open dense set specified in Lemma 6, 
then there is a  positive constant $c$ with 
$$ \|\grad \omega_n \|_{C^0} \geq c \lambda^n, $$
for all $n \in \naturals$.
}

\thm{Proof:} The first statement was proved in the paragraph
above the  theorem  so assume that
$\omega_0 \in G$ from Lemma 6 and let $S_0$ be
an essential regular strip for $\omega_0$. The case where 
$\omega_0$ has an essential regular annulus
is similar. Since $S_0$ is an essential regular strip
 there is an interval of $\omega_0$ values,
$I= [\Omega_L, \Omega_R]$,  so that for $r \in I$,
$\omega_0^{-1} (r) \cap S_0 := K_{0r}$ is an essential regular arc
and $\mu(K_{0r}) = r$. 

Now fix $n > 0$. By 
Theorem 2,  $S_n := \phi_n(S_0) $ is an essential regular strip for $\omega_n$
which is made up of the essential regular arcs $K_{nr} := \phi_n(K_{0r})$. 
If $\ell_n$ is the minimum Euclidean length of any of the $K_{nr}$,
it follows from (10) that there is a constant $k$ depending
only on the homotopy class of $S_0$ so that 
$$\ell_n \geq k \lambda^n.\eqno{(11)}$$

Since $S_n$ is a regular strip, it is equipped with a pair
of orthogonal foliations by arcs, namely, the regular
arcs $K_{nr}$ and trajectories of the gradient
flow of $\omega_n$ restricted to $S_n$. This implies that we 
may find smooth function
$a_1, a_2 : I\raw \reals$ and a diffeomorphism 
$$F: 
\{ (a, r) : r \in I\ \hbox{\rm and} 
\ a_1(r) \leq a \leq a_2(r) \} \raw S_n$$
so that:

\item{(1)} $F$ restricted to $[a_1(r), a_2(r)] \times \{r\}$ parameterizes
$K_{nr}$.

\item{(2)} For $i = 1, 2$, $F(a_i(r), r)$ 
parameterizes the intersection 
of $S_n$ with a boundary component of $M_n$.

\item{(3)} $F(\{a\} \times I)$ is a trajectory of the gradient flow of 
$\omega_n$ restricted to $S_n$.
\medskip
For such an $F$, $\frac{\partial F}{\partial a}$ is orthogonal
to $\frac{\partial F}{\partial r}$,  and since $F^{-1}$ restricted
to a $F(\{a\} \times I)$ is exactly $\omega_n$, 
$$\|\frac{\partial F}{\partial r}( a, r)\| 
= \frac{1}{\|\grad w_n (F(a, r))\|}.$$
Using the fact that $\phi_n$ preserves area
$$\eqalign{
\area(S_0) &= \area(S_n) \cr
&= 
\int_{\Omega_L}^{\Omega_R} \int_{a_1(r)}^{a_2(r)} \det(DF) \; da\; dr \cr
&= \int_{\Omega_L}^{\Omega_R} \int_{a_1(r)}^{a_2(r)} 
\frac{\|\frac{\partial F}{\partial a}\|}{ \|\grad w_t (F(a, r))\|}\; da\; dr \cr
&\geq \frac{1}{\| \grad\omega_n(z_n)\|} 
\int_{\Omega_L}^{\Omega_R} \int_{a_1(r)}^{a_2(r)} 
\|\frac{\partial F}{\partial a}\| \; da\; dr \cr
&\geq \frac{ (\Omega_R - \Omega_L) \ell_n}{\| \grad\omega_n(z_n)\|},\cr}$$
where $z_n $ is such that 
$\| \grad\omega_n(z_n)\| 
= \max \{ \|\grad \omega_n(z)\|  : z \in S_n\}$.
Thus using (11), 
$$\| \grad\omega_n(z_n)\| \geq 
\frac{ (\Omega_R - \Omega_L) k\lambda^n }
{\area(S_0)}
$$
and note that none of the constants depend on $n$.
\QED

Although we restricted attention here to the case of three
stirrers, there is a similar result for pA protocols with
more stirrers and for protocols for which the TN representative
in the isotopy class has at least one pA component.

In one sense Theorem 7 says that for pA protocols any compatible
Euler velocity fields are diverging as $t \raw\infty$ because 
$\| X \|_{C^2} \raw \infty$, or alternatively, $\|\laplacian X\|_{C^0}
\raw \infty$. Examining
the situation more carefully one sees that this divergence is a result of
the ``piling up'' of levels sets of $\omega$ and thus the 
graph of the vorticity is acquiring sharp ridges
packed closely together. Thus the vorticity is becoming
more evenly distributed and regular in the sense that 
in many places the local mean vorticities are 
approaching the global mean.
Now if the fluid motion were time-periodic
and strong mixing with
respect to Lebesgue measure, this would be the behavior of  $\mu\circ\phit$
for any $L^2$-function $\mu$. While the fluid motion is certainly not periodic
and strong mixing it is worth noting that the proof of Theorem 7 indicates that 
this behavior is the result of each time one advance map of the fluid being 
isotopic to a pA map and so is, in a certain sense, the ``memory'' of
the strong mixing of the pA map.

\section{6}{Discussion and Questions}
We begin by summarizing the contributions  of this paper 
to Question 1, and more generally, to the understanding of 
the dynamics of Euler diffeomorphisms.
The results fall into various cases based on the 
initial vorticity and the isotopy class induced by the stirrer motion. For 
expositional simplicity we continue to restrict the 
discussion to fluid regions with the topology of $M_0$ from
\S 5, but similar conclusions can be drawn for 
systems with more than $3$ stirrers.

As a starting point, recall that 
the set of area preserving diffeomorphisms with positive topological
entropy on a genus zero surface is dense and open 
in the $C^\infty$ topology ([W], see [F] for other generic
conditions).  Thus in addressing 
Question 2 we shall primarily focus
on the question of the existence of Euler diffeomorphisms with
positive entropy.

The simplest case is that of stationary boundaries, i.e. no stirring.
From Proposition 3, we know that in this case constant initial vorticity
yields only steady Euler fluid
motions and thus zero entropy. From Proposition 4
we know that whatever the stirring protocol, for typical initial
vorticity a time-periodic Euler fluid motion
also has zero entropy. Thus any Euler 
model for chaotic dynamics in this class comes from 
vorticity that is atypical and non-constant.

At the other extreme of stirring,  for pA protocols 
Theorem 7 says that typical initial vorticity never gives rise to 
Euler diffeomorphisms. 
So once again to find Euler models for chaotic dynamics we must
have atypical initial vorticity. Amongst these atypical systems, perhaps
the most attractive for additional study are those with constant
vorticity. By Proposition 3 we know that these systems do, in fact, 
always have periodic Euler velocity fields. In addition,  since the class 
is pA, the resulting Euler diffeomorphisms all have  positive entropy. 
Further, Proposition 3 also shows that for a fixed protocol
there is only a $4$-parameter family of such Euler diffeomorphisms.
Amongst these constant vorticity solutions the zero vorticity case is 
especially attractive because in that case the stream function is 
harmonic and the methods of complex
analysis can be brought to bear. These systems will be the subject of a 
subsequent paper.  Very interesting numerical results on these systems
are contained in [FCB].

Between the two extremes of stirring are the finite order protocol
in which the stirrers move in paths that are topologically
the same as those generated by circulating them in
order around a circle.  This case also requires further investigation. 

In the $C^\omega$-case the situation is somewhat clearer. For nonconstant
initial vorticity the entropy is always zero (Corollary 5). 
As already noted, for constant vorticity one always has
periodic compatible Euler velocity fields. For a pA protocol the resulting 
Euler diffeomorphism always have  positive entropy, but
stationary boundaries yield zero entropy.
For finite order protocols
presumably one can have both types of  behavior.

This situation is
somewhat reminiscent of Arnol'd's dichotomy for $C^\omega$-steady
3D Euler flows. In that case if the vorticity is generic, then
the Bernoulli function gives an integral of motion and forces
zero entropy. The non-generic case is when the vorticity in aligned
with the velocity field giving what is termed a Beltrami flow. 
These flows can have positive entropy and much progress has
been made by Etnyre and Ghrist based on the observation that 
Beltrami flows can be identified with the Reeb flows of contact
forms. In the analogy to 2D periodic Euler fluid motions the 
Beltrami case corresponds 
to zero vorticity, and there also one has nice additional structure,
namely, the stream function is harmonic.

While the entire picture is not yet clear, especially in  the $C^k$-case,
at this point it appears that Euler diffeomorphisms, especially ones 
with positive entropy, are rather rare. Thus it seems likely that 
Euler diffeomorphisms cannot manifest all the dynamical behavior of 
area preserving diffeomorphisms. 
However, we do not, as of yet, have a specific dynamical obstruction
that would give a negative answer to Question 2.

We close with some questions stimulated by Question 2. 
All the questions except the last assume the situation
studied here: planar, multi-connected regions with perhaps
moving boundaries.

\thm{Question 3:} {\it How common are time-periodic Euler velocity fields 
which are not steady? Give explicit examples of initial vorticity
distributions and periodic stirring protocols 
so that the compatible Euler solution is also periodic.} 

It is worth noting that the results in this paper have
only used the fact that the vorticity is passively transported
by Euler fluid motions and not the important additional information
that this function is actually the curl of the velocity fields generating the
motion. In particular, in certain cases the curl will make level sets
of the vorticity  ``curl''  up and not return to themselves periodically,
and so prevent the periodicity of the Euler solution.

\thm{Question 4:} {\it Study the dynamics of Euler diffeomorphism
arising from smooth, atypical initial vorticity such as
those with with areas of constant vorticity. 
Are there any chaotic, 
time-periodic Euler flows with stationary boundaries or which
are compatible with finite order protocols?}

It would also be very interesting to lessen the
regularity assumptions of this paper and study, for example, 
Questions 3 and 4 for systems whose vorticity is the indicator function of
a region with smooth boundary (this case is often referred
to as ``vortex patches''). 
It is known that such initial data gives rise to flow maps
that are H\"older homeomorphisms and 
the boundary curve of a vortex patch remains smooth throughout the
evolution. Thus the resulting Euler homeomorphisms
have special characteristics which could be valuable in
understanding their dynamics. In this context it is useful to note that
the Thurston-Nielsen theory is a theory about homeomorphisms, and
so, for example, if the boundary of a vortex patch is an essential
curve, then a compatible Euler velocity field cannot be periodic
under a pA stirring protocol.

It is also worth remarking that another kind of singular Euler fluid motion,
point vortices, do have periodic solutions with
positive entropy in the sense that there are {\it relative} periodic
solutions of the $3$-vortex problem on the cylinder for which
the motion of the vortices may be treated as stirrers with 
a protocol of pA type ([BAS2]).  Thus the induced velocity field
on the cylinder  is time-periodic after a space translation, and the 
time $T$-map has positive entropy. It is likely that such periodic solutions
also exist in the disk and plane without the geometric
phase. Can these examples be smoothed while maintaining the periodicity?
Note that Proposition 4 and the property of pA classes given in (10) 
force strong restrictions on how the smoothing can be done.

\thm{Question 5:} {\it Study Question 1 for the time $T$-diffeomorphisms of
time-periodic quasistationary Stokes flow in multi-connected region
with periodic stirring protocols.}

The quasistationary Stokes equation is derived from assumptions
at the other extreme from  those behind the Euler equation.
Instead of assuming the viscosity is negligible,  one assumes that the 
vorticity dominates and that the acceleration
or inertial term can be neglected. The velocity field is thus
required to satisfy:
$$\eqalign{   
\Delta X &=  - \grad p_t \hbox{} \cr 
\divv(\bX) &= 0 \cr
\bX  & =  \dalpha_i\ \  \hbox{\rm on the boundary.}\cr}
$$
where $p_t$ is again the pressure. Note that the boundary
conditions correspond to the fluid sticking to the boundary without slipping;
this is a consequence of the nonzero viscosity. Solutions to the 
quasistationary Stokes equation are characterized by 
the stream function $\Psi$ being
biharmonic at each time, $\Delta^2 \Psi = 0$. Thus in terms of mathematical
structure, quasistationary Stokes flows are  similar to time-periodic
zero vorticity Euler 
fluid motions.  Both are time-periodic Hamiltonian systems whose 
Hamiltonian has additional structure.
In the zero vorticity Euler case the Hamiltonian is harmonic, 
while for Stokes flow it is biharmonic. Note, however,  
the nontrivial difference in the boundary conditions. Stokes flows 
with a pA stirring protocol are studied in [BAS1], [FCB], and [V].

\thm{Question 6:} {\it In three dimensions the Helmholtz-Kelvin theorem
says that the vorticity (now a vector field) is transported.
Thus with generic initial vorticity a 3D time-periodic Euler fluid motion 
preserves a nontrivial vector field. What restrictions does
this place on the dynamics of 3D Euler diffeomorphisms?}

\thm{Acknowledgments:} The author would like to thank Hassan Aref for
pointing out that versions of Proposition 4 and Corollary 5 are
contained in [BS] and Igor Mezi\'c for the remark after Theorem 7
about the asymptotics of the vorticity.
\vfill\eject

\bigskip
\centerline{\bf Bibliography}
\medskip

\item{[AMR]} Abraham, R., Marsden, J., and Ratiu, T., {\it
Manifolds, tensor analysis, and applications}, Springer-Verlag, 
1988.

\bibitem{AK} 
  Arnol'd, V. \& Khesin, B., {\it Topological Methods in Hydrodynamics},
Springer Verlag, 1998.

\bibitem{Ba} 
  Batchelor, G.K., 
  {\it An Introduction to Fluid Mechanics},
Cambridge University Press, 1967.

\bibitem{Bi} Birman, J, {\it Braids, Links and Mapping Class Groups},
Annals of Mathematics Studies, Princeton University Press,  1975.

\item{[BL]} Birman, J. and Libgober, A., ed.,
Proceedings of the AMS-IMS-SIAM Joint Summer Research Conference
on Artin's Braid Group, {\it  Contemp.  Math.}, {\bf 78}, 1988.

\item{[BW]} Birman, J. and Williams, R., Knotted periodic orbits in
dynamical systems II: knot holders for fibered knots, {\it Contemp. Math.},
{\bf 20}, 1983, 1--60.

\item{[Bd1]} Boyland, P., Topological methods in surface dynamics,
 {\it Topology and its Applications}, {\bf 58}, 223--298, 1994.

\item{[Bd2]}  Boyland, P.,
 Fluid Mechanics and Mathematical Structures, {\it 
 An Introduction to the Geometry and Topology of Fluid Flows}, 
   (ed. R.L. Ricca), NATO-ASI Series: Mathematics, Kluwer, 105--134, 2001.

\item{[BAS1]} Boyland, P., Aref, H. and Stremler, M., Topological
fluid mechanics of stirring, {\it J. Fluid Mech.},
{\bf 403}, 277--304, 2000.

\item{[BAS2]} P. Boyland, H. Aref and M. Stremler, 
   Topological fluid dynamics of point vortex motions, {\it Physica D}, 
{\bf 175},  69--95, 2003.

 \item{[Br]} Brenier, Y.,  Topics on hydrodynamics and volume preserving maps,
 {\it Handbook of mathematical fluid dynamics II},  North-Holland, 55--86, 
 2003. 

\item{[BS]} Brown, M. \& Samelson, R., Particle motion in 
vorticity- conserving, two-dimensional incompressible flows,
{\it Phys. Fluids}, {\bf 6}, 2875--2876, 1994.

\item{[CB]} Casson, A. and Bleiler, S., {\it Automorphisms of Surfaces
after Nielsen and Thurston}, London Math. Soc. Stud. Texts, {\bf 9},
Cambridge University Press, 1988.

\item{[C]} Chemin, J.-Y., Fluides parfaits incompressibles, 
{\it Ast\'erisque},  {\bf 230}, 1995.

\bibitem{EG} 
Etnyre, J. \& Ghrist, R.,  Stratified integrals and unknots in inviscid 
flows, {\it Contemporary Mathematics}, {\bf 246}, 99-112, 1999. 

\item{[FLP]} Fathi, A, Lauderbach, F. and Poenaru, V., Travaux
de Thurston sur les surfaces, {\it  Ast\'erique}, {\bf 66-67},  1979.

\item{[FCB]}
 Finn, M., Cox, S. \& Byrne, H., 
Topological chaos in inviscid and viscous mixers, 
{\it J. Fluid Mech.}, {\bf 493}, 345--361, 2003.

\item{[F]} Franks, J., Rotation numbers and instability sets,
{\it   Bull. Amer. Math. Soc.}, {\bf  40}, 263--279, 2003.

\item{[FH]} Franks, J. and Handel, M., 
Periodic Points of Hamiltonian Surface Diffeomorphisms, preprint, 
ArXiv, math.DS/0303296.

\item{[G]} Gallavotti, G., {\it Foundations of fluid mechanics},
Springer-Verlag, 2002.

\bibitem{GK} 
  Ghrist, R. \& Komendarczyk, R., Topological features of inviscid flows, 
{\it
 An Introduction to the Geometry and Topology of Fluid Flows},
   (ed. R.L. Ricca), NATO-ASI Series: Mathematics, Kluwer, 183--202, 2001.
  
\item{[HH]}
He, C. and  Hsiao, L., 
Two-dimensional Euler equations in a time dependent domain,
{\it J. Differential Equations}, {\bf  163}, 265--291,  2000.

\item{[H]} Hirsch, M., 
 {\it Differential Topology},   Springer Verlag, 1976.

\item{[K1]} Katok, A., Bernoulli diffeomorphisms on surfaces, {\it Ann.
 Math.}, {\bf 110}, 529--547, 1979.

\item{[K2]} Katok, A., Lyapunov exponents, entropy
and periodic orbits for diffeomorphisms, {\it Publ Math IHES}, {\bf
51}, 137--173, 1980.

\item{[Kh]} Koch, H., Transport and instability for perfect fluids,
{\it Math. Ann.}, {\bf 323}, 491--523, 2002.

\item{[Ko]} 
Kozono, Hideo, On existence and uniqueness of a global classical 
solution of the two-dimensional Euler equation in a time-dependent domain,
{\it J. Differential Equations}, {\bf  57},   275--302, 1985.   

\item{[MB]} Majda, A. and Bertozzi, A., {\it Vorticity and incompressible
flow}, Cambridge University Press, 2002.

\item{[MP]} Marchioro, C. and Pulvirenti, M., {\it Mathematical
theory of incompressible non-viscous fluids}, Springer-Verlag, 1994.

\item{[Ma]} Matsumoto, Y., {\it An introduction to Morse theory},
    American Mathematical Society, 2002.

\item{[MS]} McDuff, D. \& Salamon, D., {\it Introduction to
Symplectic Topology}, Oxford University Press, 1995.

\bibitem{Mi1}
  Milnor, J., A note on curvature and the fundamental group, 
{\it J. of Diff. Geom.}, {\bf 2},  1--70, 1968.
  
\bibitem{Mi2}
  Milnor, J., 
  {\it Morse Theory},  Princeton University Press, 1969.

\bibitem{Mo} 
  Moser, J.,  
  On the volume elements on a manifold,
  {\it Trans. Amer. Math. Soc.}, {\bf 120}, 286--294, 1965.

\item{[P]} Polterovich, L., {\it The Geometry of the Group of
Symplectic Diffeomorphisms}, \hfill\break Birkh\"auser, 2001.

\item{[Sh]} Shnirelman, A., Diffeomorphisms, braids, and flows,
{\it An Introduction to the Geometry and Topology of Fluid Flows}, 
   (ed. R.L. Ricca), NATO-ASI Series: Mathematics, Kluwer, 253--270, 2001.

\bibitem{Se}
  Serrin, J.,  
  Mathematical Principles of Classical Fluid Mechanics,
  {\it Handbuch der\hfill\break  Physik}, {\bf VIII/1}, 1959.

\item{[T]} Thurston, W., On the geometry and dynamics of diffeomorphisms of
surfaces, {\it  Bull. A.M.S.}, {\bf  19}, 417--431,  1988.

\item{[V]} Vikhansky, A., Simulation of topological chaos in laminar flow,
 {\it Chaos}, 2004 (in press).

\item{[W]} Weiss, H., 
Genericity of symplectic diffeomorphisms of $S\sp 2$ with positive topological 
entropy. A remark on: "Planar homoclinic points" 
by D. Pixton and "On the generic existence of homoclinic points" 
by F. Oliveira, 
{\it J. Statist. Phys.}, {\bf  80},  1995, 481--485.

\item{[Y]} Young, L.-S., Entropy of continuous flows on compact $2$-manifolds,
{\it Topology}, {\bf 16},   469--471, 1977.

\bye